\documentclass{amsart}
\pdfoutput=1
\usepackage[utf8]{inputenc}
\usepackage{amsmath,calligra,mathrsfs}
\usepackage{amssymb}
\usepackage{amsthm}
\usepackage{tabularx}
\usepackage{hyperref}
\usepackage{tikz}
\usetikzlibrary{cd}
\everymath{\displaystyle}
\usepackage[hmargin = 1in,vmargin=1in]{geometry}

\newcommand{\ic}{\text{ic}}

\newcommand{\Ab}{\mathrm{Ab}}

\newcommand{\et}{_{\text{ét}}}

\newcommand{\Sp}{\mathrm{Sp}}

\newcommand{\tr}{\mathrm{tr}}

\newcommand{\Sh}{\mathrm{Sh}}
\newcommand{\Het}{\mathrm{H}_{\mathrm{\acute{e}t}}}

\newcommand{\Hom}{\mathrm{Hom}}

\newcommand{\ord}{\mathrm{ord}}
\newcommand{\ZZ}{\mathbb{Z}}
\newcommand{\Cl}{\mathcal{C}\ell}
\newcommand{\F}{\mathbb{F}}
\newcommand{\Oh}{\mathcal{O}}
\newcommand{\Gm}{\mathbb{G}_m}
\newcommand{\Pn}{\mathbb{P}}

\newcommand{\Q}{\mathbb{Q}}
\newcommand{\An}{\mathbb{A}}

\newcommand{\G}{\mathcal{G}}

\newcommand{\res}{\mathrm{res}}
\newcommand{\Id}{\mathrm{Id}}

\newcommand{\A}{\mathrm{A}}
\newcommand{\intHom}{\textbf{Hom}}

\newcommand{\Ho}{\mathrm{H}}

\newcommand{\br}{\mathrm{Br}}

\usepackage{scalerel}
\usepackage{stackengine,wasysym}

\newcommand{\bu}{\bullet}

\newcommand{\RHom}{\mathrm{R}\mathbf{Hom}}

\newcommand{\cl}{\mathrm{cl}}

\newcommand{\X}{\mathcal{X}}

\newcommand{\GG}{\mathbb{G}}
\newcommand{\inv}{\mathrm{inv}}

\newcommand{\codim}{\mathrm{codim}}
\usepackage[hmargin = 1in,vmargin=1in]{geometry}
\usepackage[backend=biber, sorting=anyvt]{biblatex}
\renewbibmacro{in:}{%
  \ifentrytype{article}{}{\printtext{\bibstring{in}\intitlepunct}}}

\newtheorem{theorem}{Theorem}[section]
\newtheorem{corollary}[theorem]{Corollary}
\newcommand{\Rhom}{\mathrm{RHom}}
\newtheorem{lemma}[theorem]{Lemma}
\newtheorem{prop}[theorem]{Proposition}
\theoremstyle{definition}
\newtheorem{definition}[theorem]{Definition}
\newtheorem{construction}[theorem]{Construction}

\newtheorem{example}[theorem]{Example}
\newtheorem{question}{Question}

\theoremstyle{remark}
\newtheorem{remark}[theorem]{Remark}

\newtheorem{notation}[theorem]{Notation}

\usepackage{blindtext}
\usepackage{subfiles}
\usepackage{comment}
\usepackage{lipsum}

\title{Evaluating the tame Brauer group of open varieties over local fields}
\author{victor de vries}
\date{January 2025}

\begin{document}

\maketitle

\begin{abstract}
In this document we let $U$ be a smooth variety of pure dimension $d$ over a local field $k_v$ with unit ball $\Oh_v$ and residue field $\F$ of characteristic $p>0$ and we set $n$ to be a positive integer such that $p\nmid n$. For various $u\in U(k_v)$ we study the evaluation map $u^*:\Het^2(U,\mu_n)\to \Het^2(k_v,\mu_n)$. We suppose that $U$ embeds as an open subscheme in a regular scheme $\X$ that is of finite type over $\Oh_v$. We assume that $Z:=\X\setminus U$ is a divisor and we endow it with its reduced scheme structure. We show that for $u_1,u_2\in U(k_v)$ that lift to $x_1,x_2\in \X(\Oh_v)$ we obtain the same evaluation map $u_1^*=u_2^*$ under the two conditions that first, there is an equality of reductions $\overline{x_1}=\overline{x_2}$ in $\X(\F)$ and second, that $\cl(x_1\cap Z)=\cl(x_2\cap Z)$ holds in $\Ho^{2d}_{|\overline{x}|}(Z,\mu_n^{\otimes d})$.
\end{abstract}

\section*{Introduction}
We fix a non archimedean local field $k_v$ with unit ball $\Oh_v$ and residue field $\F$ of cardinality $q$ and characteristic $p$. We set $U$ to be a smooth variety over $k_v$ that is equidimensional of dimension $d$. We define for any scheme $X$ the cohomological Brauer group of $X$, denoted $\br(X)$, by $\br(X):=\Het^2(X,\Gm)$. For $u\in U(k_v)$ the evaluation maps $s^*:\br(U)\to \br(k_v)$ are widely studied. Both the source and target for this arrow are torsion groups and the target group, $\br(k_v)$, comes equipped with a canonical isomorphism $\inv:\br(k_v)\to \Q/\ZZ$. 

An example in which one might study such evaluation maps is when $U$ has a model $U_k$ over a numberfield~$k$, such that $k_v$ is a completion of $k$ (at the place $v$). In this case understanding when two distinct points $x_1,x_2\in U(k_v)$ induce different evaluation maps $x_1^*,x_2^*:\br(U_k)\to \br(k_v)$ is of interest for studying the Brauer-Manin obstruction to strong approximation on $U$.

A known result in the study of evaluation maps is obtained in the setting that $U$ admits a regular model~$\mathcal{U}$ over $\Oh_v$ such that the special fiber, $\mathcal{U}_\F$, has a nonsingular locus $\mathcal{U}_\F^{\text{ns}}$ consisting of only one connected component. This result works for any point $u\in U(k_v)$ that lifts to a point $x\in \mathcal{U}(\Oh_v)$, whose reduction modulo $v$ is $\overline{x}$. In this situation there is the commuting diagram in which $(p')$ means $p$-primary torsion (for a proof see \cite[Proposition~5.1]{Badred}):
\begin{equation}\label{modelcase}
\begin{tikzcd}
\br(U)(p') \arrow[d, "u^*"] \arrow[r] & {\Ho^1(\mathcal{U}_\F^{\text{ns}},\Q/\ZZ)(p') \arrow[d, "\overline{x}^*"]} \\
\br(k_v)(p')   \arrow[r, "\sim"]                & {\Ho^1(\F,\Q/\ZZ)(p')}            
\end{tikzcd}
\end{equation}

The aim of this document is to investigate what happens when many of the conditions above are removed. For example, we may want to evaluate a point $u\in U(k_v)$ that does not lift to an $\Oh_v$-point on the model~$\mathcal{U}$.\par 
The setting in which we will be working is as follows: We let $U$ be a smooth variety over $k_v$ and we fix a regular scheme $\X$ that is of finite type over~$\Oh_v$ together with an open immersion $U\to \X$ such that the complement of $U$ is a divisor. Let $Z$ be this complement endowed with its reduced scheme structure. Any point $x\in \X(\Oh_v)$ that comes from a point $u\in U(k_v)$ has that the scheme-theoretic intersection $Z\cap x$ is an $\Oh_v/(\pi_v^m)$-point on $Z$ for certain $m\geq 1$. A weak variant of the theorem that we prove is as follows:

\begin{theorem}\label{slogan}
Let $u\in U(k_v)$ be a point that lifts to $\X(\Oh_v)$. The evaluation map $\br(U)(p')\to \br(k_v)(p')$ depends only on the closed subscheme $Z\cap x$ of $Z$.
\end{theorem}

By commuting Diagram (\ref{modelcase}) we know that this theorem holds in the case that $\X$ is a model of $U$ over~$\Oh_v$ together with condition on the special fiber in the paragraph before Diagram (\ref{modelcase}).\par
The hope that the evaluation map $u^*$ should only depend on the reduction $\overline{x}$ as in the setting of Diagram~(\ref{modelcase}) is easily crushed by taking for $p\neq 2$, $U=\GG_{m,\Q_p}$ and $\X=\mathbb{A}^1_{\ZZ_p}=\Sp(\ZZ_p[t])$. Namely for $u\in \ZZ_p^*\setminus (\ZZ_p^*)^2$ the points $p,up\in \GG_{m,\Q_p}(\Q_p)$ both reduce to $0\in \mathbb{A}^1_{\ZZ_p}(\F_p)$ and evaluating quaternion algebras $(t,a)_{-1}\in \br(\GG_{m,\Q_p})[2]$ for various $a\in \Q_p^*/(\Q_p^*)^2$ on the different points $p,up$ produces different results (see example~\ref{babyexample}).

We wish to state our main theorem, which is a stronger variant of Theorem \ref{slogan}. Denote by $|Y|$ the underlying topological space of a scheme $Y$. The closed subscheme $Z\cap x$ of $Z$ has a cohomology class $\cl(Z\cap x)\in \Ho^{2d}_{|Z\cap x|}(Z,\mu_n^{\otimes d})$. This cohomology class appears in our main theorem as follows:

\begin{theorem}\label{MainTHIntro}
The evaluation map $u^*:\br(U)[n]\to \br(k_v)[n]$ depends only on $\cl(Z\cap x)\in \Ho^{2d}_{|Z\cap x|}(Z,\mu_n^{\otimes d})$.
\end{theorem}

The intuition that the author would like to give is that the evaluation of a point $u$ corresponding to $x$ only depends on how $Z\cap x$ looks `locally' on $Z$. The above theorem gives such a statement, since in the case that $Z\cap x$ is nonempty, the excision property in étale cohomology implies that we may replace $\Ho^{2d}_{|Z\cap x|}(Z,\mu_n^{\otimes d})$ by $\Ho^{2d}_{|Z\cap x|}(A,\mu_n^{\otimes d})$, where $A$ is the henselisation of $\Oh_{Z,|Z\cap x|}$ at its maximal ideal.\\

Our initial motivation for proving such a result came from the examples that were encountered in Section~\ref{1}. Our initial idea to use the leray spectral sequence associated to $U\hookrightarrow \X$ to prove Theorem \ref{slogan}, because the $E_2$-terms in this spectral sequence are related to the cohomology of $Z$. It turns out that this method does not give the desired result and we explain why at the end of the Section \ref{1}. However if one assumes that $Z$ is an snc divisor and that $\gcd(n,q-1)=1$ holds, then the spectral sequence argument does work as is also explained at the end of section \ref{1}.

Our proof of Theorem \ref{MainTHIntro} is inspired by the methods used in Section 5 of \cite{FailureWA}. Like in this reference we will use the language of derived categories. To make the document more accessible to people that might not be so familiar with derived categories, we give an overview of the necessary tools from derived categories in Section \ref{2}.
Section \ref{3} is dedicated to the proof of Theorem \ref{MainTHIntro}.

We note that studying the evaluation maps on the $p^\infty$ torsion $u^*:\br(U)\{p\}\to \br(k_v)\{p\}$ is drastically harder than studying the evaluation maps on the prime to $p$ torsion. For the case in which $\X=\mathcal{U}$ is a smooth model of $U$ with geometrically connected special fiber, we refer to the paper by M. Bright and R. Newton (\cite{WildBr}) for the study of the evaluation maps $u^*:\br(U)\to \br(k_v)$ for points $u\in U(k_v)$ that lift to~$\mathcal{U}(\Oh_v)$ (so also the $p^\infty$-torsion is evaluated here). 

\section{Concrete examples and the insufficiency of the Leray spectral sequence}\label{1}
In this section we exhibit a class of examples to Theorem~\ref{slogan}. Later in the section we explain why a certain spectral sequence argument fails to prove Theorem \ref{slogan} in general, but we show that the argument goes through under additional hypothesis. The notation of this section will be as in Theorem~\ref{slogan}. We start with a `baby-example' to which we already referred in the introduction.

\begin{example}\label{babyexample}
Set $\X=\An^1/\ZZ_p$ and fix the open subset $U=\Gm/\Q_p$ of $\X$. Since we have $\br(\An^1_{\Q_p})=\br(\Q_p)$, the purity theorem \cite[Theorem 3.7.1]{BrGroGroBook} gives an exact sequence 
\begin{equation}
\label{BrGm}0\to \br(\Q_p)(p')\to \br(U)(p')\to \Ho^1(\Q_p,\Q/\ZZ)(p')\to 0
\end{equation}
(the exact sequence ends in $0$ because of the vanishing of $\Het^3(\An^1_{\Q_p},\Gm)$).
For any $n,m\in \ZZ$ with $p\nmid n$, $a\in \Q_p^\times$ and $b\in \Ho^1(\Q_p,\ZZ/n)$, we can consider the element $[a\cdot X^m]_n\cup b\in \br(U)[n]$ where with $[-]_n$ we mean its class in $\Het^1(K(U),\mu_n)$.
By combining formula (1.27) in \cite{BrGroGroBook} for the residue map $\br(k(U))\to \Ho^1(\Q_p,\Q/\ZZ)$ with exact sequence \ref{BrGm}, we obtain that $\br(U)(p')$ is generated by elements of the form $[a\cdot X^m]_n\cup b\in \br(U)[n]$.\par
So two points $u_1,u_2\in U(\Q_p)$ yield the same evaluation map $\br(U)(p')\to \br(\Q_p)(p')$ if and only if for all $n$ with $p\nmid n$ and for all $b\in \Ho^1(\Q_p,\ZZ/n)$ the elements $[u_1]_n\cup b$ and $[u_2]_n\cup b$ agree.
By Poitou-Tate duality (\cite[p.28]{ArithmeticDuality}) this is equivalent to for all $n$ with $p\nmid n$ there being an equality $[u_1]_n=[u_2]_n$ in $\Ho^1(\Q_p,\mu_n)=\Q_p^\times/(\Q_p^\times)^n$. Writing $u_1=v_1\cdot p^{l_1}$ and $u_2=v_2\cdot p^{l_2}$ with $v_i\in \ZZ_p^\times$ gives that the condition $[u_1]_n=[u_2]_n$ for all $n>0$ is equivalent to us having both $l_1=l_2$ and $v_1\equiv v_2\,(\text{mod }(1+p))$ in $\ZZ_p^\times$.\par
Denote by $Z$ the divisor $\Sp(\ZZ_p[X]/(pX))\subset \X$. The complement of $Z$ is $U$. Assume that $x_1,x_2$ lift to $x_1,x_2\in \X(\ZZ_p)$ or equivalently that we have $l_1,l_2\geq 0$. The conditions $l_1=l_2$ and $v_1\equiv v_2\,(\text{mod }(1+p))$ are equivalent to there being an equality of ideals $(pX,X-u_1)=(pX,X-u_2)$, so we see that $u_1$ and $u_2$ induce the same evaluation map $\br(U)(p')\to \br(\Q_p)(p')$ if and only if we have $Z\cap x_1=Z\cap x_2$.
\end{example}

We also have the following example in which we evaluate `obvious elements' in the Brauer group of a variety. 

\begin{example}
Let $U$ be a variety over $k_v$, which embeds as an open subset in a variety $\X$ over $\Oh_v$. Let $f,g$ be elements of $\Oh_{\X}(\X)$ that become invertible in $\Oh_{\X}(U)$. Denote by abuse of notation for any scheme~$Y$ and for any positive integer $n$ by $N_n$ the norm map $\br(Y(\zeta_n))\to \br(Y)$. For any $p\nmid n$ the cyclic algebra $N_n(\{f,g\}_{\zeta_n})$ is an element of $\br(U)$. Evaluating such an element at a point $u\in U(k_v)$ gives us the cyclic algebra $N_n(\{f(u),g(u)\}_{\zeta_n})$ which is an element of $\br(k_v)$. Write $f(u)=v_f\cdot \pi_v^{m_f}$ and write $g(u)=v_g\cdot \pi_v^{m_g}$ for $v_f,v_g\in \Oh_v^\times$. Under the residue isomorphism $\br(k_v)\to \Ho^1(\F_v,\Q/\ZZ)$ the element $N_n(\{f(u),g(u)\}_{\zeta_n})$ is mapped to $N_{\F_v(\zeta_n)/\F_v}[(-1)^{m_fm_g}\frac{v_f^{m_g}}{v_g^{m_f}}]_{\zeta_n}$, where with $[-]$ we mean the class in $\Ho^1(\F_v(\zeta_n),\mu_n)$.\par
Assume that $u$ lifts to $x\in \X(\Oh_v)$ and so we have $m_f,m_g\geq 0$. Denote for any $h\in \Oh_{\X}(\X)$ by $Z_h$ the closed subscheme $V(\pi_v\cdot h)$ of $\X$.
Like in the previous example we have that the expression $v_f\cdot \pi_v^{m_f}$ (resp. $v_g\cdot \pi_v^{m_g}$) is determined by $Z_f\cap x$ (resp. $Z_g\cap x$) up to multiplying $v_f$ (resp. $v_g$) by an element in $(1+\pi_v)\Oh_v^\times$.
So therefore the evaluation of $N_n(\{f,g\}_{\zeta_n})$ under a point $u\in U(k_v)$ corresponding to $x\in \X(\Oh_v)$ depends only on the intersection $x\cap Z_f$ in $Z_f$ and on the intersection $x\cap Z_g$ in $Z_g$.
This in turn is determined by $x\cap Z_{fg}$ inside $Z_{fg}$.
In particular, for $V:=\X\setminus Z_{fg}$ and for $B$ the subgroup of $\br(V)$ spanned by elements of the form $N_n\{f^i,g^j\}_{\zeta_n}$, we have gained insight in the pairing $B\times V(k_v)\to \br(k_v)$.
\end{example}

\subsection*{The Leray spectral sequence}
We let $U,\X,Z,u,x$ be given as before. They sit in the following cartesian diagram
\begin{equation} \label{LSSCD}
\begin{tikzcd}
U \arrow[r, "j"]                        & \X                        \\
\Sp(k_v) \arrow[u, "u"] \arrow[r, "j'"] & \Sp(\Oh_v) \arrow[u, "x"]
\end{tikzcd}   . 
\end{equation}
We obtain Leray spectral sequences associated to $j$ and $j'$, which are $E_2^{s,q}:=\Het^s(\X,R^qj_*\mu_n)\implies \Het^{s+q}(U,\mu_n)$ and $A_2^{s,q}:=\Het^s(\Oh_v,R^qj'_*\mu_n)\implies \Het^{s+q}(k_v,\mu_n)$. Moreover, the commuting Diagram (\ref{LSSCD}) provides us with a morphism of spectral sequences $E_r^{s,q}\to A_r^{s,q}$. Our hope is now as follows: We wish to understand the map $u^*:\Het^2(U,\mu_n)\to \Het^2(k_v,\mu_n)$ in terms of the maps $E_2^{p,q}\to A_2^{p,q}$. Namely the terms $E_2^{p,q}$ should be related to the cohomology of $Z$ on which the $\Oh_v/(\pi_v^n)$-point $x_Z:=x\cap Z$ lies. So our hope is that understanding the pullback map $x_Z^*:\Het^*(Z)\to \Het^*(\Oh_v/(\pi_v^n))$ gives us full understanding of the pullback map $u^*$.\\

For the spectral sequence corresponding to $j'$ we get a filtration $\Het^2(k_v,\mu_n)=:F_0\supset F_1\supset F_2\supset F_3:=\{0\}$. One easily deduces that one has $F_2=\br(\Oh_v)[n]=0$ and $F_0=F_1=A^{1,1}_2\cong \Het^1(\F_v,\ZZ/n)$, where the isomorphism of the last group with $\Het^2(k_v,\mu_n)$ is a very special case of the following proposition due to Gabber:

\begin{prop}\label{SNC}
Consider a regular $\ZZ[1/n]$-scheme $X$ with a strict normal crossings divisor $D$, whose components are $\iota_i:D_i\to X$ and whose complement is $j:U\to X$. There is a canonical isomorphism $\res:R^1j_*\ZZ/n\to \bigoplus_{i}\iota_{i,*}(\ZZ/n(-1))$ and $R^*j_*\ZZ/n$ is the exterior algebra on $R^1j_*\ZZ/n$.
\end{prop}

Suppose for the moment that the complement $j:U\to \X$ is an snc divisor $\cup_i D_i$ as in the proposition. There is also a filtration $H^2(U,\mu_n)=:V_0\supset V_1\supset V_2\supset V_3:=\{0\}$ induced by the spectral sequence, where we have $V_0/V_1=E_{\infty}^{0,2}$. By the proposition we have $E_2^{0,2}=\bigoplus_{i<j}\ZZ/n(-1)(D_i\cap D_j)$, which is nonzero as soon there are at least two $D_i$ whose intersection $D_i\cap D_j$ is nonempty and $\gcd(\#\F-1,n)\neq 1$ holds. The following example shows that even $E_{\infty}^{0,2}=0$ might be nonzero, so we may have $\Het^2(U,\mu_n)=V_0\supsetneq V_1$.

\begin{example}\label{FiltrationIncompatible}
Take $U$ and $\X$ as in Example \ref{babyexample}, denote by $j:U\to \X$ be the natural map and let $n$ be such that $\gcd(n,p-1)\neq 1$. The Leray spectral sequence has $E_2$ terms $E_2^{3,0}=\Het^3(\An^1_{\ZZ_p},\mu_n)=0$ and $E_2^{2,1}=\Het^2(\An^1_{\F_p},\ZZ/n)\oplus \Het^2(\Sp(\ZZ_p),\ZZ/n)=0$, where for calculating $E_2^{2,1}$ we used Proposition \ref{SNC}. The same proposition gives us that we have $E_2^{2,0}=\Hom(\mu_n(\F_p),\ZZ/n)$, which is nonzero. The vanishing of $E_2^{2,1}$ and $E_2^{3,0}$ shows that we have $E_2^{0,2}=E_{\infty}^{0,2}$ and so we have $\Het^2(U,\mu_n)=:V_0\supsetneq V_1$.
\end{example}

The problem that one has whenever $\Het^2(U,\mu_n)=:V_0\supsetneq V_1$ holds (which in general can not be avoided by Example \ref{FiltrationIncompatible}) is that one can not understand the evaluation map $u^*:\Het^2(U,\mu_n)\to \Het^2(k_v,\mu_n)$ by investigating the induced maps in the spectral sequences $E_{\infty}^{2-q,q}\to A_{\infty}^{2-q,q}$. Indeed, these induced maps are the maps $V_i/V_{i+1}\to F_i/F_{i+1}$ from which one can not deduce what $u^*\alpha\in \Het^2(k_v,\mu_n)$ is for some $\alpha\in V_0\setminus V_1$. Therefore we will work in the derived category, which will allow us to not `lose information by taking homology too quickly'.

\begin{remark}\label{when E_2^{0,2}=0}
It should be noted that in some natural examples we do have $\Het^2(U,\mu_n)=:V_0=V_1$ and thus the evaluation map $u^*:\br(U)[n]\to\br(k_v)[n]$ is completely determined by the map between the terms on the $E_2$-pages in degree $(1,1)$, which is $\Het^1(\X,R^1j_*\mu_n)\to \Het^1(\Oh_v,R^1j'_*\mu_n)=\Het^1(\F_v,\ZZ/n)$. One of these cases was already mentioned, namely when the complement of $U$ is an snc divisor and when we have $(n,\#\F-1)=1$, since then $E_2^{0,2}$ vanishes. We will prove in the theorem below that in this case the evaluation map $u^*$ can be understood by looking at residues. 

Another case in which $E_2^{0,2}$ vanishes is when the complement of $U$ is an snc divisor, whose components do not meet. In this case we can change $\X$ by removing all the components of the complement of $U$ that do not meet the lift of $u$ to $\X$, so that the complement of $U$ equals $\X_\F$, which is integral. Then we can apply Diagram (\ref{modelcase}) to understand the evaluation map $u^*$.
\end{remark}

Now we state the theorem that was mentioned in the last remark.

\begin{theorem}\label{snc theorem} Assume that $Z:=\X\setminus U$ is an snc divisor with components $(Z_i)_i$. For $n$ a positive integer such that $\gcd(q-1,n)=\gcd(n,p)=1$ the following holds.
\begin{itemize}
\item The image of the residue map $\Het^2(U,\mu_n)\overset{(\res_i)_i}{\to} \bigoplus_{i}\Het^1(K(Z_i),\ZZ/n)$ is contained in the subgroup $\bigoplus_{i}\Het^1(Z_i,\ZZ/n)$ of the target. \item Denote by $n_i$ the length of the scheme $B_i:=Z_i\times_{\X,x}\Sp(\Oh_v)$ and denote whenever $B_i$ is nonempty by $x_i:\Sp(\F)\hookrightarrow Z_i$ the inclusion of the reduced subscheme of $B_i$ into $Z_i$.  The following diagram commutes:
\begin{equation}
\begin{tikzcd}
{\Het^2(U,\mu_n)} \arrow[d, "u^*"] \arrow[r, "(\res_i)_i"] & {\bigoplus_i \Het^1(Z_i,\ZZ/n)} \arrow[d, "\sum_i n_i x_i^*"] \\
{\Het^2(k_v,\mu_n)} \arrow[r, "\res"]                & {\Het^1(\F,\ZZ/n)}                 
\end{tikzcd}
\end{equation}
\end{itemize}
\begin{proof}As mentioned in the last remark, the assumptions imply that $E_2^{0,2}$ vanishes. We have a diagram in which the left square commutes since this square is induced by the morphism of Leray spectral sequences:
\begin{equation}\label{leraycd}
\begin{tikzcd}
{\Het^2(U,\mu_n)/F_0} \arrow[d, "u^*"] & {\Het^1(\X,R^1j_*\mu_n)} \arrow[l, "\sim"'] \arrow[d, "x^*"] \arrow[r, "\sim"] & {\bigoplus_i \Het^1(Z_i,\ZZ/n)} \arrow[d, "\sum_i n_ix_i^*"] \\
{\Het^2(k_v,\mu_n)}                    & {\Het^1(\Oh_v,R^1j'_*\mu_n)} \arrow[l, "\sim"'] \arrow[r, "\sim"]             & {\Het^1(\F_v,\ZZ/n)}               
\end{tikzcd}
\end{equation}

The residue map $\Het^2(U,\mu_n)\overset{(\res_i)_i}{\to} \bigoplus_{i}\Het^1(K(Z_i),\ZZ/n)$ can be constructed as follows (\cite{BrGroGroBook}[Theorem 3.7.1]): First shrink $\X$ and $Z$ by removing all the intersections $Z_i\cap Z_j$ with $i\neq j$. Then apply Proposition \ref{SNC} to the new $\X,Z$. So by naturality of the residue map we obtain that $\Het^2(U,\mu_n)\overset{(\res_i)_i}{\to} \bigoplus_{i}\Het^1(K(Z_i),\ZZ/n)$ extends $\Het^2(U,\mu_n)\overset{(\res_i)_i}{\to} \bigoplus_{i}\Het^1(Z_i,\ZZ/n)$.

To prove the second bullet point in the theorem it remains to be proven that the right square of Diagram~(\ref{leraycd}) commutes. For this it will suffice to show that following diagram of morphisms of sheaves commutes:
\begin{equation}
\begin{tikzcd}\label{sheafsdiag}
R^1j_*\mu_n \arrow[d, "x^*"] \arrow[r, "\res"] & \bigoplus_i (\iota_i)_*\ZZ/n \arrow[d, "\sum_in_ix_i^*"] \\
x_*R^1j'_*\mu_n \arrow[r, "\res"]              & (x\circ \iota_\F)_*\ZZ/n                                
\end{tikzcd}
\end{equation}

Here the rightmost vertical arrow is on the $i$'th component, evaluated on an étale $W\to \X$, defined by the map $\ZZ/n^{\pi_0(W\times_\X D_i)}\to \ZZ/n^{\pi_0(\Sp(\F)\times_\X W)}$ that maps the standard basis vector $e_B$ to $(n_i\cdot \delta_{B,C})_C$, where~$C$ ranges over $\pi_0(\Sp(\F)\times_\X W)$, $B$ ranges over $\pi_0(W\times_\X D_i)$ and $\delta_{B,C}$ is $1$ if $C\subset B$ and $0$ else. The stalk of the bottom-right term is only nonzero at a geometric point $\tilde{x}$ of $\X$ lying over $x\circ \iota_\F$. Note that for any morphism $f:Y\to S$, the sheaf $R^mf_*\mathcal{F}$ is the sheafification of the presheaf $\Het^m(-\times_SY,\mathcal{F})$. So to prove that Diagram (\ref{sheafsdiag}) commutes it will suffice to show that the following diagram commutes:
\begin{equation}
\begin{tikzcd}\label{passtolimitdiag}
{\varinjlim\Het^1(W\times_\X U,\mu_n)} \arrow[d, "u^*"] \arrow[r, "\res"] & \varinjlim\bigoplus_i \ZZ/n^{\pi_0(W\times_\X Z_i)} \arrow[d, "\sum_in_ix_i^*"] \\
{\varinjlim\Het^1(W\times_\X \Sp(k_v),\mu_n)} \arrow[r, "\res"]            & \varinjlim\ZZ/n^{(\pi_0(W\times_\X \Sp(\F))}                     
\end{tikzcd}
\end{equation}
Here the colimit is taken over connected étale $(W,\tilde{w})\to (\X,\tilde{x})$. Given any such $(W,\tilde{w})\to (\X,\tilde{x})$, we can first pass to the limit over all Zariski refinements of $W$ (so taking $\Sp(\Oh_{W,w})$, where $\tilde{w}$ lies over a closed point $w$ of $W$) and then take étale covers of $\Sp(\Oh_{W,w})$ to obtain the same colimit. Since $W$ is connected and regular (because $\X$ is regular), the local ring $\Oh_{W,w}$ is a UFD. Let $\{D_{i_j}\}$ be the components of $Z$ in $\Sp(\Oh_{W,w})$ (where we use $i_j$ to indicate that $D_{i_j}$ lies over $Z_i$) and write $A$ for the subring of $\Oh_{W,w}$ corresponding to the complement of $\cup_{i_j} D_{i_j}$. Then we have $\Het^1(A,\mu_n)=A^\times/(A^\times)^n$ by the Kummer sequence. Denote by $K$ the field $k_v\otimes_{\Oh_v}\Oh_{W,w}$ and by $\Oh$ its ring of integers. The field $K$ is an unramified extension of $k_v$ and there is an are induced points $x':\Sp(\Oh)\to \Oh_{W,w}$ and $u':\Sp(K)\to \Sp(A)$. For proving that Diagram (\ref{passtolimitdiag}) commutes it hence suffices to show that the following commutes (where $n_{i_j}:=n_i$):
\begin{equation}
\begin{tikzcd}\label{ringdiag}
A^\times/(A^\times)^n \arrow[r, "(\ord_{D_{i_j}})_{i_j}"] \arrow[d, "u'^*"] & \bigoplus_{i_j} \ZZ/n \arrow[d, "\sum_{i_j}n_{i_j}"] \\
K^\times/(K^\times)^n \arrow[r, "\ord_v"]                            & \ZZ/n                                
\end{tikzcd}
\end{equation}
Let $f_{i_j}\in A_W$ cut out to $D_{i_j}$ and consider an element $f:=\prod_{i_j}f_{i_j}^{m_{i_j}}\cdot u$ of $A^\times$, where $u$ doesn't vanish at $w$, i.e. we have $\ord_v\pi(u)\neq 0$. Going down and then right in the diagram maps $[f]$ to $\sum_{i_j} m_{i_j}\cdot \ord_v(f_{i_j}(u'))$. Going right and then down in the diagram produces $\sum_{i_j}m_{i_j}\cdot n_{i_j}$. These terms agree since we have \begin{equation}\ord_v(f_{i_j}(u'))=\mathrm{length}(D_{i_j}\times_{\Oh_{W,w},x'} \Sp(\Oh))=\mathrm{length}(Z_i\times_{\X,x}\Sp(\Oh))=:n_i=:n_{i_j}\end{equation} Here we used in the second equality that $W\to \X$ is étale. We conclude that Diagram (\ref{ringdiag}) commutes. This was the remaining step for proving the theorem.
\end{proof}
\end{theorem}

We go back to our `baby-example' to exhibit an element of $\Het^2(U,\mu_n)$ that does not come from the first step in the filtration, $V_1$.

\begin{remark}
One can remove the assumption on $n$ in the theorem if one replaces $\Het^2(U,\mu_n)$ by $V_1$, the first step in its Leray filtration. The proof of this is identical.

Now we go back to Example \ref{FiltrationIncompatible}, where we take $p$ to be odd and we take $n=2$. We consider $\alpha\in \Het^2(U,\mu_2)$ that maps to the quaternion algebra $(p,t)_{-1}\in \br(U)[2]$. Consider the $\Q_p$-points $u_1:=p$ and $u_2:=-p$ on~$U$. They lift to $\ZZ_p$-points on $\X$ and both points intersect $V(p)$ and $V(t)$ with multiplicity $1$ and moreover $p$ and~$-p$ both reduce to $0$ modulo $p$. So if we would have $\alpha\in V_1$, then Theorem \ref{snc theorem} would tell us that we have $u_1^*(\alpha)=u_2^*(\alpha)$. However, this is not the case, since $u_i^*(\alpha)$ can be identified with $u_i^*(p,t)_{-1}\in \br(\Q_p)$ and we have $u_1^*(p,t)_{-1}=(p,p)_{-1}$, which corresponds under the residue isomorphism to $-1\in \F_p^*/(\F_p^*)^2$, while $u_2^*(p,t)_{-1}=(p,-p)_{-1}$ corresponds to $1\in \F_p^*/(\F_p^*)^2$.
\end{remark}

Note that actually the values of this residue isomorphism should be in $\ZZ/2$, but since we already picked a primitive second root of unity there is no ambiguity.

\section{Derived categories}\label{2}
Let $X$ be a scheme and set $n$ to be a positive integer which is invertible on $X$. Denote by $\Lambda$ the constant sheaf $\ZZ/n$ on $X\et$. In this section we will work with the category $D^+(X\et,\Lambda)$, which is also called the \textit{bounded above derived category of} $\Lambda$ \textit{modules} on~$X\et$. Recall that this is a \textit{triangulated category}. For a more thorough introduction on derived categories in general we refer the reader either chapter 13 of the Stacks Project (\cite[\href{https://stacks.math.columbia.edu/tag/05QI}{Tag 05QI}]{stacks-project}) or any textbook on this subject such as \cite{BookDerivedCategories}. There are three subsections, the first is used for introducing the \textit{localisation triangle}, the second for introducing \textit{Gysin maps} corresponding to local complete intersection morphisms (see also \cite{TravauxdeGabber}[Exposé XVI.2]) and the last for introducing \textit{representable functors}.

\subsection{The localisation triangle}\label{localisationtriangle}
We now wish to introduce the \textit{localisation triangle}. For this let $\iota:Z\to X$ be a closed immersion with open complement $j:U\to X$. We have the functors on the étale sites of sheaves $j^*,j_*,\iota_*\iota^!$ of which $j^*$ and $\iota_*$ are exact. The canonical adjoint pairs $(j^*,j_*)$ and $(\iota_*,\iota^!)$ lift to adjoint pairs $(j^*,Rj_*)$ and $(\iota_*,R\iota^!)$ on the derived categories (see \cite[\href{https://stacks.math.columbia.edu/tag/0FNC}{Tag 0FNC}]{stacks-project}). In particular we obtain natural transformations $\eta_j:\Id\to Rj_*j^*$ and $\epsilon_\iota:\iota_*R\iota^!\to \Id$ between endofunctors on $D^+(X\et,\Lambda)$. Now we give the localization triangle.

\begin{lemma}\label{localisation}
The above adjunctions induce for any complex $K^\bu$ a distinguished triangle \begin{equation}\label{Triangle1}\iota_*R\iota^!K^\bu\to K^\bu\to Rj_*j^*K^\bu\to \iota_*R\iota^!K^\bu[1].\end{equation}
\begin{proof}
Take a quasi-isomorphism $K^\bu\to\mathcal{I}^\bu$ with $\mathcal{I}^\bu$ a complex of injectives. By using the definition of $j_*j^*$ and by applying (\cite[\href{https://stacks.math.columbia.edu/tag/01EA}{Tag 01EA}]{stacks-project}) we see that for each $i$ the map $\mathcal{I}^i\to j_*j^*\mathcal{I}^i$ is surjective. Thus we have an exact sequence of injective complexes $0\to \iota_*\iota^!\mathcal{I}^\bu\to \mathcal{I}^\bu\to j_*j^*\mathcal{I}^\bu\to 0$, which is termwise split because on the level of sheaves the functors $\iota_*,\iota^!,j_*,j^*$ preserve injectives. And so the distinguished triangle is obtained.
\end{proof}
\end{lemma}

We have the following lemma that gives a criterion for a morphism of triangles to exist.

\begin{lemma}\label{uniquetrianglemorphism}(\cite{PerverseSheaves}[Proposition 1.1.9])
For a diagram as follows, where the rows are distinguished triangles and the middle square commutes, the following are equivalent:
\begin{equation*}
\begin{tikzcd}
K^\bu \arrow[r, "f"] \arrow[d, dotted, "\alpha"] & L^\bu \arrow[r, "g"] \arrow[d, "\beta"] & M^\bu \arrow[r, "h"] \arrow[d, "\gamma"] & {K^\bu[1]} \arrow[d, dotted, "\alpha{[1]}"] \\
K'^\bu \arrow[r, "f'"]                 & L'^\bu \arrow[r, "g'"]                  & M'^\bu \arrow[r, "h'"]                   & {K'^\bu[1]}                 
\end{tikzcd}
\end{equation*}
\begin{itemize}
    \item It holds that $g'\circ \beta \circ f=0$.
    \item There exists a morphism $\alpha:K^\bu\to K'^\bu$ such that the dotted arrows give a morphism of triangles.
\end{itemize}
If we have moreover that $\Hom_{D^+}(K^\bu,M'^\bu[-1])=0$, the hypothetical arrow $\alpha$ is uniquely determined by the property that it makes the leftmost square commute.
\end{lemma}

An important thing for us is how the localisation triangle behaves functorially. Namely let $U'\overset{j'}{\hookrightarrow}X'\overset{\iota'}{\hookleftarrow} Z'$ to be another triple as in the paragraph below Subsection \ref{localisationtriangle}. Given a morphism $f:X'\to X$ such that $f$ restricts to $f_U:U'\to U$ and $f_Z:Z'\to Z$. Let $K^\bu$ be an object of $D^+(X\et,\Lambda)$. By taking the localisation triangle on $X'$ of $f^*K^\bu$ and applying $Rf_*$ to this triangle, we get by a distinguished triangle on $D^+(X\et,\Lambda)$: 
\begin{equation}\label{Triangle2}
Rf_*\iota'_*R\iota'^!f^*K^\bu\to Rf_*f^*K^\bu\to R(f\circ j')_*(f\circ j')^*K^\bu\to Rf_*\iota'_*R\iota'^!f^*K^\bu[1].
\end{equation}

The following lemma describes the functoriality between triangle (\ref{Triangle1}) and triangle (\ref{Triangle2}).

\begin{lemma}\label{morphismlocalisation}
There is a unique morphism from triangle (\ref{Triangle1}) to triangle (\ref{Triangle2}) such that $K^\bu\to f_*f^*K^\bu$ is the morphism induced by adjunction $\eta_f:\Id\to Rf_*f^*$ and such that $Rj_*j^*K^\bu\to R(f\circ j')_*(f\circ j')^*K^\bu$ is the morphism induced by $Rj_*\eta_{f_U}j^*:Rj_*j^*\to  R(f\circ j')_*(f\circ j')^*$, where $\eta_{f_U}:\Id\to R(f_{U})_*f_{U}^*$ is the adjunction morphism.
\begin{proof}
The middle diagram in the hypothetical morphism of triangles is induced by the following diagram of natural transformations depicted below. 
\begin{equation}
\begin{tikzcd}
\Id \arrow[r, "\eta_j"] \arrow[d, "\eta_f"] & Rj_*j^* \arrow[d, "Rj_*\eta_{f_U}j^*"] \\
Rf_*f^* \arrow[r, "Rf_*\eta_{j'}f^*"]       & R(f\circ j')_*(f\circ j')^*           
\end{tikzcd}
\end{equation}
This diagram commutes since both ways around it equals the adjunction morphism $\Id\to R(f\circ j')_*(f\circ j')^*$ (see \cite{Categoriesforworkingmathematician} p.103).
We may replace $K^\bu$ and $f^*K^\bu$ by a complexes of injectives $\mathcal{I}^\bu$ and $\mathcal{J}^\bu$ respectively. We have by (\cite[\href{https://stacks.math.columbia.edu/tag/05TG}{Tag 05TG}]{stacks-project}) the equality $\Hom_{D^+}(\iota_*R\iota^!K^\bu,R(f\circ j')_*(f\circ j')^*K^\bu)=\Hom_{K^+}(\iota_*\iota^!\mathcal{I}^\bu,j_*(f_U)_*j'^*\mathcal{J}^\bu)$. The same statement is true if we shift the target complex by $[-1]$ in both $\Hom$'s. By Lemma \ref{uniquetrianglemorphism} it suffices to show that the previous $\Hom$-group vanishes to get a morphism of triangles and by the same lemma we have that the morphism of triangles is unique if the same $\Hom$-group vanishes but with $\mathcal{J}^\bu$ replaced by $\mathcal{J}^\bu[-1]$.
So it suffices to show that for any sheaves $\mathcal{F}\in (Z\et,\Lambda)$ and $\mathcal{G}\in (U\et,\Lambda)$ we have $\Hom(\iota_*\mathcal{F},j_*\mathcal{G})=0$. This is clear by applying the adjunction $(j^*,j_*)$ because $j^*\iota_*\mathcal{F}$ is the zero-sheaf.
\end{proof}
\end{lemma}

\subsection{Gysin maps}\label{SecGysin}
All morphisms in this section are assumed to be of finite type and all schemes are assumed to be Noetherian. In this subsection we introduce the \textit{Gysin maps} corresponding to local complete intersection morphisms of $\ZZ[1/n]$-schemes with an ample invertible sheaf (see \cite[\href{https://stacks.math.columbia.edu/tag/01PS}{Tag 01PS}]{stacks-project} for the definition of the last term). For more details we refer to \cite{TravauxdeGabber}[Exposé XVI.2].

We have the following definition, which is well defined (see \cite[\href{https://stacks.math.columbia.edu/tag/0F7I}{Tag 0F7I}]{stacks-project}).
\begin{definition}
Let $f:X\to Y$ be a morphism of schemes. Let $j:X\to \overline{X}$ be an open immersion with a proper morphism $\overline{f}:\overline{X}\to Y$ extending $f$. Define $Rf_!:D^+(X\et,\Lambda)\to D^+(Y\et,\Lambda)$ to be the functor $R\overline{f}_*\circ j_!$ where $j_!$ is the (exact) extension by zero functor.
\end{definition}

In particular we see that $Rf_!=Rf_*$ if $f$ is proper.

\begin{remark}
It is known that $Rf_!$ admits a right adjoint $Rf^!:D^+(Y\et,\Lambda)\to D^+(X\et,\Lambda)$ (see \cite[\href{https://stacks.math.columbia.edu/tag/0G2C}{Tag 0G2C}]{stacks-project}). The notation $Rf^!$ might be a bit misleading since in general $Rf^!$ is \textbf{not} a derived functor of some functor $(Y\et,\Lambda)\to(X\et,\Lambda)$. This is however the case if $f=\iota$ is a closed immersion, so that $\iota^!:(Y\et,\Lambda)\to(X\et,\Lambda)$ is defined and or in the case that $f=j$ is étale, in which case we have $Rj^!=j^*$. 
\end{remark}

We make a remark that says that Tate twists commute with almost all functors in this paper.

\begin{remark}
Associated to a morphism $f:X\to Y$ we have seen the set of functors $S=\{Rf_*,f^*,Rf_!,Rf^!\}$ in this paper. All of these functors commute with taking Tate twists, i.e. for any object $K^\bu$ (in the appropriate derived category of modules over $\Lambda$) and for any
$F\in S$ we have $F(K^\bu\otimes \Lambda(d))=F(K^\bu)\otimes \Lambda(d)$. We will use this throughout the paper.\end{remark}

We now work towards defining the Gysin morphism.

\begin{definition}
A \textit{regular immersion} $\iota:Z\to X$ is an immersion such that for $U\subset X$ an open such that $Z\to U$ is a closed immersion, the ideal sheaf $\mathcal{I}_Z\subset \Oh_U$ is locally given by a regular sequence. 
\end{definition}

Now we define a local complete intersection morphism.

\begin{definition}\label{lci}
A \textit{local complete intersection morphism} (abr. \textit{lci morphism}) is a morphism $f:X\to Y$ such that every $x\in X$ has an open neighbourhood $U$ such that $f|_U$ admits a factorization $f|_U:U\overset{\iota}{\to} P\overset{q}{\to}Y$ where~$\iota$ is a regular immersion and $q$ is smooth.
Define $\mathcal{I}^{\text{ic}}$ to be the category whose objects are $\ZZ[1/n]$-schemes that admit an ample invertible sheaf and whose morphisms are lci morphisms.
\end{definition}

Of course the definition of $\mathcal{I}^{\ic}$ requires that the composition of two lci morphisms is an lci morphism, which is true (see \cite[\href{https://stacks.math.columbia.edu/tag/069J}{Tag 069J}]{stacks-project}) using `independence' of chosen local factorizations as in the remark below).

\begin{remark}\label{rklci}
Suppose that we know that around $x$ there is a local factorization $f|_U:U\overset{\iota}{\to} P\overset{q}{\to}Y$ where $\iota$ is a regular immersion and $q$ is smooth. Then locally around $x$ any factorization $f|_V:V\overset{\iota'}{\to} P' \overset{q'}{\to} Y$ with $q'$ smooth and $\iota'$ an immersion yields that $\iota'$ is a regular immersion (\cite[\href{https://stacks.math.columbia.edu/tag/069E}{Tag 069E}]{stacks-project}). So in particular if an lci morphism $X\to Y$ factors globally $X\overset{\iota}{\to} P\overset{q}{\to} Y$, with $q$ smooth and $\iota$ an immersion, then $\iota$ is automatically a regular immersion. If the lci morphism is a morphism in $\mathcal{I}^{\ic}$, then $X$ and $Y$ having ample invertible sheaves implies such a global factorization $X\to \Pn^m_Y\to Y$.
\end{remark}

By doing a local computation, using the above remark and (\cite[\href{https://stacks.math.columbia.edu/tag/00NR}{Tag 00NR}]{stacks-project}), it follows that a morphism $X\to Y$ between regular schemes is an lci morphism.

\begin{definition}
To a morphism $f:X\to Y$ in $\mathcal{I}^{\ic}$ with a global factorization $X\overset{\iota}{\to} P\overset{q}{\to} Y$ as in Remark~\ref{rklci} we define the \textit{relative virtual dimension} of $f$ to be the integer $\text{rel.}\dim(q)-\codim(\iota)$.
\end{definition}

Now we make the following definition for easier notation of the Gysin maps as in \cite[Exposé XVI]{TravauxdeGabber}.

\begin{definition}
Let $f:X\to Y$ be a morphism in $\mathcal{I}^\ic$ with relative virtual dimension  $d$. Denote by $f^{?}$ the functor $D^+(Y\et,\Lambda)\to D^+(X\et,\Lambda)$ that sends a complex $K^\bu$ to $K^\bu\mapsto Rf^!K^\bu(-d)[-2d]$.
\end{definition}

If we have morphisms $X\overset{f}{\to} Y\overset{g}{\to} Z$, all in $\mathcal{I}^\ic$, then it it is clear that we have $f^?\circ g^?=(g\circ f)^?$. We give the construction of the Gysin morphism of a morphism in $\mathcal{I}^{\ic}$. 

\begin{construction}
Let $f:X\to Y$ be a morphism in $\mathcal{I}^{\ic}$ of relative virtual dimension $c$ with a global factorization $f:X\overset{\iota}{\to}P\overset{q}{\to}Y$ as in Remark \ref{rklci} such that $\iota$ has codimension $c$. We will define a morphism in $D^+(X\et,\Lambda)$ denoted $\cl_f:\Lambda\to f^?\Lambda$. After defining this morphism for the cases $f=\iota$ and $f=q$ we can define $\cl_f$ to be the following composition: $\Lambda\overset{\cl_{\iota}}{\to} \iota^? \Lambda \overset{\iota^?\cl_q}{\to} R(q\circ \iota)^?\Lambda$,
\begin{itemize}
    \item Set $f=\iota$ to be a regular immersion of codimension $c$. Then $f$ factors as $j\circ \overline{\iota}$ where $j$ is an open immersion and $\overline{\iota}$ is a closed immersion. We have $j^!=j^*$ and so we have a canonical isomorphism $R\overline{\iota}^!\Lambda\to R\iota^!\Lambda$ induced by the isomorphism $\Lambda\to j^*\Lambda$. Therefore we may assume that $\iota$ is a closed immersion.
    
    In the case that $\iota$ is a closed immersion of codimension $c$, the semi-purity theorem tells us that $R^e\iota^!\Lambda(c)=0$ for $e<2c$ (see e.g. \cite{FailureWA}[Lemma 5.4]). This implies by (\cite[\href{https://stacks.math.columbia.edu/tag/06XS}{Tag 06XS}]{stacks-project}) that taking homology in degree $0$ gives an isomorphism of groups $\Hom_{D^+}(\Lambda,R\iota^?\Lambda)\to \Hom_{\Sh}(\Lambda,R^{2c}\iota^!\Lambda(c))$. A~morphism of sheaves $\Lambda\to R^{2c}\iota^!\Lambda$ is determined by picking a global section of $R^{2c}\iota^!\Lambda(c)$. By the mentioned semi-purity theorem, these global sections are precisely $\Ho_X^{2c}(P,\Lambda(c))$ and in here we have a distinguished element $\cl(X)$ corresponding to the cycle class of $X$ in $P$ (\cite{SGA4.5} p.140). We set $\cl_\iota:\Lambda\to \iota^!\Lambda(c)[2c]$ to be the morphism that sends $1$ to $\cl(X)$ after taking $h^0$.
    \item Set $f=q$ to be a smooth morphism of relative dimension $e$. We have $R^aq_!\Lambda(d)=0$ for $a>2e$ (by using \cite[\href{https://stacks.math.columbia.edu/tag/0DDF}{Tag 0DDF}]{stacks-project}). Therefore we get group isomorphisms depicted below of which the first is by adjunction $(Rq_!,Rq^!)$ and the second is by taking homology in degree $0$ (again using \cite[\href{https://stacks.math.columbia.edu/tag/06XS}{Tag 06XS}]{stacks-project}): \[\Hom_{D^+}(\Lambda,q^?\Lambda)\to \Hom_{D^+}(Rq_!\Lambda(e)[2e],\Lambda)\to \Hom_{\Sh}(R^{2e}q_!\Lambda(e),\Lambda)\]
    In this last group we have a distinguished element called the `trace map' $\tr_q:R^{2e}q^!\Lambda(e)\to \Lambda$ (see \cite{SGA4} Exposé XVIII Théorème 2.9). We set $\cl_q:\Lambda\to Rq^!(-e)[-2e]$ to be the morphism that corresponds to $\tr_q$ above.
\end{itemize}
\end{construction}

Of course we want the construction to be independent of the choice of $\iota$ and $q$. 

\begin{lemma}\label{compgys}{\cite{TravauxdeGabber}[Exposé XVI.2 Théorème 2.5.5 and Théorème 2.5.12]}
The construction of $\cl_f$ does not depend on the chosen factorization $f=q\circ \iota$. For $f:X\to Y$ and $g:Y\to Z$ morphisms in $\mathcal{I}^{\ic}$ of relative virtual dimension $d$ and $e$ respectively, the composition $\Lambda\overset{\cl_f}{\to} f^?\Lambda\overset{f^?\cl_g}{\to} (g\circ f)^?\Lambda$ equals~$\cl_{g\circ f}$.
\end{lemma}

Let $f:X\to Y$ be a morphism in $\mathcal{I}^\ic$. Since $n$ is invertible on $Y$ we have a canonical isomorphism $\Lambda_X\cong~f^*\Lambda_Y$. The following lemma shows that the Gysin morphism $\Lambda\to f^?\Lambda$ extends to a natural transformation.

\begin{lemma}\label{gysinnat}
Let $f:X\to Y$ be a morphism in $\mathcal{I}^\ic$ of relative virtual dimension $d$. There is a natural transformation $\Cl_f:f^*\to f^?$ such that when evaluated on $\Lambda$ we have $\Cl_f(\Lambda)=\cl_f:f^*\Lambda\to f^?\Lambda$.
\begin{proof}
We need to give for any object  $K^\bu$ of $D^+(Y\et,\Lambda)$ a map $f^*K^\bu\to f^?K^\bu$, that is natural in $K^\bu$. We have the equality $K^\bu=K^\bu\otimes^L\Lambda$ for any such object $K^\bu$. Since pullback and derived tensor product commute we have $f^*K^\bu=f^*K^\bu\otimes^L f^*\Lambda$.
We can use $\cl_f:f^*\Lambda\to f^?\Lambda$ on the second component to get a map $f^*K^\bu\otimes^L f^*\Lambda\to f^*K^\bu\otimes^Lf^?\Lambda=[f^*K^\bu \otimes^L f^!\Lambda](-d)[-2d]$.
It is clear that this map is natural in $K^\bu$.
So it suffices to give a natural transformation $f^*(-)\otimes^Lf^!\Lambda\to f^!(-)$.\par
We have the compactly supported projection formula $(-)\otimes^LRf_!(-)\overset{\sim}{\to}Rf_![f^*(-)\otimes^L(-)]$, see \cite[\href{https://stacks.math.columbia.edu/tag/0GL5}{Tag 0GL5}]{stacks-project} (that it is a natural transformation follows from that in the reference they use the usual projection formula combined with a natural transformation $(-)\otimes^Lj_!(-)\to j_![j^*(-)\otimes^L(-)]$). Inverting the compactly supported projection formula and using the adjunction $(Rf_!,f^!)$ gives us a natural transformation of bifunctors: \[f^*(-)\otimes^L(-)\to f^![(-)\otimes^LRf_!(-)]\]
So we obtain a map $f^*K\otimes^Lf^!\Lambda\to f^![K^\bu\otimes^LRf_!f^!\Lambda]$ that is natural in $K^\bu$. Postcomposing this map with the map $f^![K^\bu\otimes^LRf_!f^!\Lambda]\to f^!K^\bu$ induced by the adjunction map $\epsilon^!_{f}:Rf_!f^!\Lambda\to \Lambda$ produces the natural transformation on the component $K^\bu$.
\end{proof}
\end{lemma}

\subsection{Representable functors}
In this subsection we introduce the notion of a \textit{representable functor} on $D^+(X\et,\Lambda)$.

\begin{remark}{\textit{(Internal Hom)}}
There is a functor $\intHom(\mathcal{F},-):(X\et,\Lambda)\to (X\et,\Lambda)$ for $\mathcal{F}$ a sheaf of $\Lambda$-modules. This functor evaluated on a sheaf $\G$ gives the sheaf which in turn evaluated on an étale $f:U\to X$ gives $\Hom(f^*\mathcal{F},f^*\mathcal{G})$. This functor can be `extended' to a functor $\intHom^\bu(K^\bu,-):K^+(X\et,\Lambda)\to K^+(X\et,\Lambda)$ for any complex $K^\bu$ (i.e. when $K^\bu=\mathcal{F}$ is a single sheaf and the input is a sheaf $\mathcal{G}$ we have $\intHom^\bu(\mathcal{F},\mathcal{G})=\intHom(\mathcal{F},\mathcal{G})$ as a complex concentrated in degree $0$). This functor $\intHom(K^\bu,-)$ should be thought of as the functor that when evaluated on an object $L^\bu$ produces a complex of sheaves, whose $i$'th component classifies the homomorphisms of sheaves $K^n\to L^{n+i}$ for all $n\in \ZZ$ (see (\cite[\href{https://stacks.math.columbia.edu/tag/0A8K}{Tag 0A8K}]{stacks-project}). 

The functor $\intHom^\bu(K^\bu,-)$ in turn `lifts' to a functor $\RHom(K^\bu,-):D^+(X\et,\Lambda)\to D^+(X\et,\Lambda)$ which when evaluated on injective complexes equals $\intHom^\bu(K^\bu,-)$. In the case that $K^\bu=\mathcal{F}$ is a sheaf we have that $\RHom(\mathcal{F},-)$ is the right derived functor of $\intHom(\mathcal{F},-)$.
\end{remark}

\begin{definition}
A functor $G:D^+(X\et,\Lambda)\to D^+(X\et,\Lambda)$ is \textit{represented by} a complex $K^\bu$  if there is an isomorphism of functors $G\to \RHom(K^\bu,-)$.
\end{definition}

\begin{example}We have the following examples of representable functors.
\begin{itemize}
    \item The identity functor is represented by the sheaf $\Lambda$, which follows from the fact that for any complex~$K^\bu$ the equality $\intHom(\Lambda,K^\bu)=K^\bu$ holds.
    \item If $G:D^+(X\et,\Lambda)\to D^+(X\et,\Lambda)$ is represented by $K^\bu$, then $G\circ (m)[n]$ (Tate twist and shift) is represented by $K^\bu(-m)[-n]$.
\end{itemize}
\end{example}

We can add another derived functor to the list of examples by using the following lemma. We assume that it is known, however we could not find a reference.

We write $\Lambda_Z$ for $\iota_*\iota^*\Lambda$, where $\iota:Z\hookrightarrow X$ is a closed immersion.

\begin{lemma}\label{representHZ}
Let $\iota:Z\hookrightarrow X$ be a closed immersion with open complement $j:U\hookrightarrow X$. The functor $\iota_*R\iota^!$ is represented by $\Lambda_Z$.
\begin{proof}
Since the functor $\iota^!:(X\et,\Lambda)\to (X\et,\Lambda)$ has an exact left-adjoint $\iota_*$, it holds that $\iota^!$ preserves injectives. So the natural map $R(\iota_*\iota^!)\to \iota_*R\iota^!$ is an isomorphism. It will hence suffice to give an isomorphism $\iota_*\iota^!\to \intHom(\Lambda_Z,-)$. There is an isomorphism $\iota_*\iota^!\cong \mathcal{H}_Z$, where $\mathcal{H}_Z(\mathcal{F}):=\ker(\mathcal{F}\to j_*j^*\mathcal{F})$.

It is easy to see how the adjunction between $\mathcal{H}_Z$ and $\iota_*\iota^*$ works. Namely the kernel of $\mathcal{F}\to \iota_*\iota^*\mathcal{F}$ is precisely $j_!j^*\mathcal{F}$. Thus this gives us a commuting diagram of solid arrows with exact rows:
\begin{equation}\label{adjunctionHZ}
\begin{tikzcd}
0 \arrow[r] & j_!j^*\Lambda \arrow[r] & \Lambda \arrow[r] \arrow[d, "(1)"', dashed] & \iota_*\iota^*\Lambda \arrow[r] \arrow[d, "(2)", dashed] & 0                 \\
            & 0 \arrow[r]             & \mathcal{H}_Z(\mathcal{F}) \arrow[r]                  & \mathcal{F} \arrow[r]                                    & j_*j^*\mathcal{F}
\end{tikzcd}
\end{equation}
Since any map $j_!j^*\Lambda\to \mathcal{H}_Z(\mathcal{F})$ is zero (look at the stalks), an arrow $(1)$ in the diagram induces a unique arrow $(2)$ giving a commuting square. Conversely, since any map $\iota_*\iota^*\Lambda\to j_*j^*\mathcal{F}$ is zero, an arrow $(2)$ induces a unique arrow $(1)$ giving a commuting square.

Let $f:V\to X$ be an étale morphism. For $a:Y\to X$ write $a_f:Y\times_XV\to V$ for the base-change along~$f$. Applying the exact functor $f^*$ to the terms in Diagram (\ref{adjunctionHZ}) gives us $f^*(j_*j^*\mathcal{F})\cong (j_f)_*j_f^*f^*\mathcal{F}$ (since $f$ is smooth) and so $f^*\mathcal{H}_Z(\mathcal{F})=\mathcal{H}_{f^*(Z)}(f^*\mathcal{F})$. Moreover since $\iota$ is proper we have $f^*\iota_*\iota^*\Lambda\cong (\iota_f)_*\iota_f^*f^*\Lambda$. The sheaf $j_*j^!\Lambda$ is characterized as the largest subsheaf with support on $U$, which by an easy argument shows that $f^*j_*j^!\Lambda=(j_f)_!j_f^*f^*\Lambda$. As a result we get that applying $f^*$ to Diagram (\ref{adjunctionHZ}) yields:
\begin{equation}{\label{basechangeofadjunction}}
\begin{tikzcd}
0 \arrow[r] & (j_f)_!(j_f)^*f^*\Lambda \arrow[r] & f^*\Lambda \arrow[d, "(1)"', dashed] \arrow[r] & (\iota_f)_*\iota_f^*f^*\Lambda \arrow[d, "(2)", dashed] \arrow[r] & 0                            \\
            & 0 \arrow[r]                        & \mathcal{H}_{f^*(Z)}(f^*\mathcal{F}) \arrow[r] & f^*\mathcal{F} \arrow[r]                                          & (j_f)_*(j_f)^*f^*\mathcal{F}
\end{tikzcd}
\end{equation}
Now we show how to define an isomorphism of sheaves $\Phi_\mathcal{F}:\mathcal{H}_Z(\mathcal{F})\to \intHom(\Lambda_Z,\mathcal{F})$. We have a canonical identification $\mathcal{H}_Z(\mathcal{F})=\intHom(\Lambda,\mathcal{H}_Z(\mathcal{F}))$. For any étale $f:V\to X$ we need to give an isomorphism of groups $\Phi_{\mathcal{F}}(f):\Hom(f^*\Lambda,f^*\mathcal{H}_Z(\mathcal{F}))\to \Hom(f^*\Lambda_Z,f^*\mathcal{F})$. Take for this isomorphism the one obtained from Diagram~(\ref{basechangeofadjunction}). Given an étale morphism $f':V'\to X$ such that $f=f'\circ g$ for some $g:V\to V'$. We need to show that $\Phi_\mathcal{F}(f)$ and $\Phi_\mathcal{F}(f')$ are compatible with the restriction map $g^*$. This follows because both the solid diagrams corresponding to $f$ (resp. $f'$) as in Diagram (\ref{basechangeofadjunction}) are obtained from applying $f^*$ (resp. $(f')^*$) to Diagram (\ref{adjunctionHZ}). So applying $g^*$ gives the required compatibility.

Last we need to show that the components $\Phi_\mathcal{F}$ give $\Phi:\intHom(\Lambda,\mathcal{H}_Z(-))\to \intHom(\Lambda_Z,-)$, a morphism of functors. This is because we used an adjunction to define the morphism, which is natural in the second argument.
\end{proof}
\end{lemma}

Now we give an important lemma for us which may be seen as a `Derived' Yoneda lemma.

\begin{lemma}
Let $F,G:D^+(X\et,\Lambda)\to D^+(X\et,\Lambda)$ be represented by $K^\bu,L^\bu$ respectively. There is a canonical isomorphism of groups $\Hom(F,G)\to \Hom_{D^+}(L^\bu,K^\bu)$.
\begin{proof}
We will give maps in both directions that are inverse to each other. We have by (\cite[\href{https://stacks.math.columbia.edu/tag/0B6E}{Tag 0B6E}]{stacks-project}) that for any two complexes~$A^\bu$ and $B^\bu$ the equality $R\Gamma(X, \RHom(A^\bu,B^\bu))=\Rhom(A^\bu,B^\bu)$ holds. Applying $h^0$ to this last term gives $\Hom_{D^+}(A^\bu,B^\bu)$. Let $\Phi:\RHom(K^\bu,-)\to \RHom(L^\bu,-)$ be a homomorphism of functors. Plugging in $K^\bu$ and applying $h^0\circ R\Gamma(X,-)$ yields a homomorphism of groups $\Hom_{D^+}(K^\bu,K^\bu)\to \Hom_{D^+}(K^\bu,L^\bu)$. The image of the identity $\Id_{K^{\bu}}$ thus defines a morphism $\alpha(\Phi):L^\bu\to K^\bu$ in the derived category. We assign to $\Phi$ this element $\alpha(\Phi)$.

For the map in the other direction: Let $f:L^\bu\to K^\bu$ be a morphism in $D^+(X\et,\Lambda)$. Let $\mathcal{I}^\bu$ be an injective complex. We need to give a map $\RHom(K^\bu,\mathcal{I}^\bu)\to\RHom(L^\bu,\mathcal{I}^\bu)$. We have by (\cite[\href{https://stacks.math.columbia.edu/tag/08J7}{Tag 08J7}]{stacks-project}) that $\RHom(K^\bu,\mathcal{I}^\bu)=\intHom^\bu(K^\bu,\mathcal{I}^\bu)$ and that $\RHom(L^\bu,\mathcal{I}^\bu)=\intHom^\bu(L^\bu,\mathcal{I}^\bu)$. Representing $K^\bu$ by a complex of injectives $\mathcal{J}^\bu$ then gives that $\intHom^\bu(K^\bu,\mathcal{I})\cong \intHom^\bu(\mathcal{J}^\bu,\mathcal{I})$ (\cite[\href{https://stacks.math.columbia.edu/tag/0A95}{Tag 0A95}]{stacks-project}). The morphism $L^\bu\to \mathcal{J}^\bu$ is represented by a morphism of complexes (in a homotopy class). So we get a morphism $\intHom^\bu(L^\bu,\mathcal{I}^\bu)\to \intHom^\bu(\mathcal{J}^\bu,\mathcal{I}^\bu)$ and so a morphism $\RHom(K^\bu,\mathcal{I}^\bu)\to\RHom(L^\bu,\mathcal{I}^\bu)$ which we call $\beta(f)$.
A tedious manipulation of definitions shows that $\alpha$ and $\beta$ are inverses of each other, proving the lemma.
\end{proof}
\end{lemma}

\section{Proof of the main theorem} \label{3}
In this section we wish to prove the main theorem of this document, which is Theorem \ref{MainTH}. We set some notation for this section. 

\begin{notation}
We let $U$ be a smooth variety over $k_v$ with structure map $q'$. We assume that there is a regular scheme $\X$ over with a finite type morphism $q:\X\to \Sp(\Oh_v)$ and an open immersion $j:U\hookrightarrow \X$, whose complement is a divisor. Denote by $Z$ this complement with its reduced scheme structure and denote the inclusion by $\iota:Z\hookrightarrow \X$. Fix a point $x\in \X(\Oh_v)$ such that the generic point $u$ of $x$ lies in $U(k_v)$. Let $x_Z:B\to Z$ denote the base-change of $x$ along $\iota$. This is summarized in the following commuting diagram of which the two leftmost squares are cartesian:
\begin{equation} \label{setting}
\begin{tikzcd}
B \arrow[d, "\overline{\iota}"] \arrow[r, "x_Z"'] & Z \arrow[d, "\iota"] &            \\
\Sp(\Oh_v) \arrow[r, "x"]                         & \X \arrow[r, "q"]    & \Sp(\Oh_v)  \\
\Sp(k_v) \arrow[u, "j'"] \arrow[r, "u"]           & U \arrow[u, "j"]  \arrow[r, "q'"]    &       \Sp(k_v) \arrow[u, "j'"]   
\end{tikzcd}
\end{equation}
\end{notation}

\begin{remark}\label{lci remark}
The assumptions on $\X$ imply that $x$ and $q$ are morphisms in $\mathcal{I}^\ic$ as in Definition \ref{lci} (see Remark \ref{rklci}). Moreover since $x$ is an lci morphism, so is the flat base-change $u$ (equivalently $u$ is a smooth point on $U$). Since $Z$ is a divisor and since $\X$ is regular, the ring $\Oh_{Z,y}$ is a complete intersection local ring, and in particular thus Cohen-Macaulay, for any $y\in Z$. From this and from the fact that $B$ is either zero-dimensional or empty, it follows that $x_Z$ is an lci morphism. So for $f$ any of the above morphisms in Diagram (\ref{setting}) it makes sense to write $\cl_f:f^*\Lambda\to f^?\Lambda$ and $\Cl_f:f^*\to f^?$ as introduced in Subsection \ref{SecGysin}.
\end{remark}

\subsection{Setup for the morphism of triangles}
One can consider the localisation triangle (see Lemma \ref{localisation}) attached to the sheaf $q^*\mu_n$ and $(Z,\X,U)$, which is as follows. We denote this triangle by $\mathcal{T}_\X$ and it is shown below:
\begin{equation}
\begin{tikzcd}
\iota_*\iota^!q^*\mu_n \arrow[r, "\epsilon_{\iota}^!"] & q^*\mu_n \arrow[r, "\eta_j"] & Rj_*j^*q^*\mu_n \arrow[r] & {\iota_*\iota^!q^*\mu_n[1]}
\end{tikzcd}
\end{equation}

Moreover we can consider the localisation triangle attached to $\mu_n$ and $(B,\Sp(\Oh_v),\Sp(k_v))$. Then we apply~$x_*$ to this triangle to get a triangle $\mathcal{T}_{\Oh_v}$, which is shown below:
\begin{equation}
\begin{tikzcd}
x_*\overline{\iota}_*\overline{\iota}^!\mu_n \arrow[r, "x_*\epsilon_{\overline{\iota}}^!"] & x_*\mu_n \arrow[r, "x_*\eta_{j'}"] & x_*Rj'_*(j')^*\mu_n \arrow[r] & {x_*\overline{\iota}_*\overline{\iota}^!\mu_n[1]}
\end{tikzcd}
\end{equation}

By applying Lemma \ref{morphismlocalisation}, we find a unique morphism of triangles $\mathcal{T}_\X\to \mathcal{T}_{\Oh_v}$ with the following properties (notation as in Lemma \ref{uniquetrianglemorphism}):
\begin{itemize}
    \item The map $\beta:q^*\mu_n\to x_*\mu_n$ is induced by adjunction $\eta_x:\Id\to x_*x^*$
    \item The map $\gamma:Rj_*j^*q^*\mu_n\to x_*Rj'_*(j')^*\mu_n$ is induced by the adjunction morphism $\eta_{u}:\Id\to u_*u^*$. More precisely, it is the morphism $Rj_*\eta_u$ applied to $Rj_*j^*q^*\mu_n$: 
    \[Rj_*j^*q^*\mu_n\to Rj_*u_*u^*j^*q^*\mu_n=x_*Rj'_*(j')^*\mu_n\]
\end{itemize}

We denote by $\alpha$ the morphism $\iota_*\iota^!q^*\mu_n\to x_*\overline{\iota}_*\overline{\iota}^!\mu_n$ in the morphism $\mathcal{T}_\X\to\mathcal{T}_{\Oh_v}$ above. The whole point of considering this `enlargement' $\X$ of $U$ is the following:

\begin{lemma}\label{determinedbyalpha}
The evaluation map $u^*:\Ho^2(U,\mu_n)\to \Ho^2(k_v,\mu_n)$ is determined by the component $\alpha$ of the morphism $\mathcal{T}_\X\to \mathcal{T}_{\Oh_v}$ (still using the notation from Lemma \ref{uniquetrianglemorphism}).
\begin{proof}
The evaluation map $x^*$ is obtained by applying $R\Gamma(\X,-)$ to the map $\gamma:Rj_*j^*q^*\mu_n\to x_*Rj'_*(j')^*\mu_n$ and then taking second homology. Applying $R\Gamma(\X,-)$ to the morphism of triangles $\mathcal{T}_\X\to \mathcal{T}_{\Oh_v}$ produces a new morphism of triangles in $D^+(\Ab)$. Taking the long exact sequences in cohomology of these new triangles gives a commuting diagram as follows:
\begin{equation}
\begin{tikzcd}
... \arrow[r] & {\Het^2(\X,\mu_n)} \arrow[r, "j^*"] \arrow[d, "x^*"] & {\Het^2(U,\mu_n)} \arrow[r, "\delta"] \arrow[d, "u^*"] & {\Ho^3_{|Z|}(\X,\mu_n)} \arrow[r] \arrow[d] & ...\\
... \arrow[r] & {\Het^2(\Oh_v,\mu_n)} \arrow[r, "(j')^*"]               & {\Het^2(k_v,\mu_n)} \arrow[r, "\delta'"]                & {\Ho^3_{|B|}(\Oh_v,\mu_n)} \arrow[r]      & ...
\end{tikzcd}
\end{equation}
It suffices to show that $\delta':\Het^2(k_v,\mu_n)\to \Ho^3_{|B|}(\Oh_v,\mu_n)$ is an isomorphism. Indeed if this is true, then we have $u^*=(\delta')^{-1}\circ h^3(R\Gamma(\X,-)(\alpha))\circ \delta$, which is determined by $\alpha$. Since $\Ho^i_{|B|}(\Oh_v,-)$ is the $i$'th right derived functor of $\mathcal{F}\mapsto \ker(\mathcal{F}(\Oh_v)\to \mathcal{F}(k_v))$, we see that we have $\Ho^i_{|B|}(\Oh_v,-)=\Ho^i_{|\Sp(\F)|}(\Oh_v,-)$. The Gysin isomorphism (\cite[p.242]{TravauxdeGabber}) $\Ho^i(\F,\ZZ/n)\to \Ho^{2+i}_{|\Sp(\F)|}(\Oh_v,\mu_n)$ allows us to replace $\Ho^3_B(\Oh,\mu_n)$ by $\Ho^1(\F,\ZZ/n)$. The map $\Ho^2(k_v,\mu_n)\to \Ho^1(\F_v,\ZZ/n)$ is the residue isomorphism, of which $\br(k_v)\overset{\sim}{\to}\Ho^1(\F_v,\Q/\ZZ)$ is perhaps a more familiar form.
So we conclude that $\Ho^2(k_v,\mu_n)\to \Ho^3_{|B|}(\Oh_v,\mu_n)$ is an isomorphism and thus that the lemma is true.
\end{proof}
\end{lemma}

\begin{remark}
All of the functors that appear in the distinguished triangles $\mathcal{T}_\X$ and $\mathcal{T}_{\Oh_v}$ commute with Tate twists. Moreover we have for any objects $K^\bu$, $L^\bu$ of $D^+(\X\et,\Lambda)$ that $\Hom(K^\bu,L^\bu)=\Hom(K^\bu(1),L^\bu(1))$ where we write `$=$' because the isomorphism is just $f\mapsto f\otimes 1$. So for understanding the map $\alpha$ in $\mathcal{T}_\X\to \mathcal{T}_{\Oh_v}$ we may replace both $\mathcal{T}_\X$ and $ \mathcal{T}_{\Oh_v}$ by the same triangles but with $\mu_n$ replaced with $\Lambda$.
\end{remark}

Denote these new triangles involving $\Lambda$ by $\Sigma_{\X}$ and $\Sigma_{\Oh_v}$. They come with a canonical morphism $\Sigma_{\X}\to \Sigma_{\Oh_v}$ of which we also denote the map $\iota_*\iota^!q^*\Lambda\to x_*\overline{\iota}_*\overline{\iota}^!\Lambda$ by $\alpha$ by an abuse of notation.

\subsection{Understanding the map $\alpha$}

We now wish to understand the map $\alpha:\iota_*\iota^!q^*\Lambda\to x_*\overline{\iota}_*\overline{\iota}^!\Lambda$ that occurs in the morphism of triangles $\Sigma_\X\to \Sigma_{\Oh_v}$. By Lemma \ref{uniquetrianglemorphism} the arrow $\alpha$ is the unique arrow such that the following diagram commutes:
\begin{equation}
\begin{tikzcd}
\iota_*\iota^!q^*\Lambda \arrow[r, "\epsilon_{\iota}^!"] \arrow[d, "\alpha"]                & q^*\Lambda \arrow[d, "\eta_{x}"] \\
x_*\overline{\iota}_*\overline{\iota}^!\Lambda \arrow[r, "\epsilon^!_{\overline{\iota}}"] & x_*\Lambda       
\end{tikzcd}
\end{equation}

First we find a seemingly more complicated way to write the arrow $\eta_x$ in this diagram. This will turn out to be useful later.

\begin{lemma}\label{chifactorization}
The above map $\Lambda_{\X}\to x_*\Lambda_{\Oh_v}$ in $D^+(\X\et,\Lambda)$ factors as in the following diagram:
\begin{equation}
\begin{tikzcd}
{\Lambda_{\X}} \arrow[r, "\cl_q"] \arrow[d, "\eta_x"] & {q^?\Lambda_{\Oh_v}} \arrow[rd, "\eta_x"]     &                                                  \\
{x_*\Lambda_{\Oh_v}}                                   & {x_*x^?q^?\Lambda_{\Oh_v}} \arrow[l, "="] & {x_*x^*q^?\Lambda_{\Oh_v}} \arrow[l, "x_*\Cl_x"]
\end{tikzcd}
\end{equation}
\begin{proof}
Since $\eta_x:\Id\to x_*x^*$ is a natural transformation, there is a commuting square:
\begin{equation}
\begin{tikzcd}
\Lambda_{\X} \arrow[d, "\eta_x(\Lambda)"] \arrow[r, "\cl_q"]   & q^?\Lambda_{\Oh_v} \arrow[d, "\eta_{x}(q^?(\Lambda))"] \\
x_*x^*\Lambda_{\X} \arrow[r, "x_*x^*\cl_q"'] & x_*x^*q^?\Lambda_{\Oh_v}      
\end{tikzcd}
\end{equation}
Therefore it is sufficient to show that $x_*\Cl_x\circ x_*x^*\cl_q:x_*\Lambda_{\Oh_v}=x_*x^*q^*\Lambda_{\Oh_v}\to x_*x^?q^?\Lambda_{\Oh_v}=x_*\Lambda_{\Oh_v}$ is the identity map. Since $\Cl_q:x^*\to x^?$ is a natural transformation, there is a commuting square:
\begin{equation}\label{chi,q,gysin}
\begin{tikzcd}
x^*q^*\Lambda_{\Oh_v} \arrow[d, "x^*\cl_q"] \arrow[r, "\cl_x"] & x^?q^*\Lambda_{\Oh_v} \arrow[d, "x^?\cl_q"] \\
x^*q^?\Lambda_{\Oh_v} \arrow[r, "\Cl_x"]           & x^?q^?\Lambda_{\Oh_v}          
\end{tikzcd}
\end{equation}
By Lemma \ref{compgys} we have that the going right and then down in Diagram (\ref{chi,q,gysin}) produces a map $\Lambda_{\Oh_v}\to \Lambda_{\Oh_v}$ equal to the identity. So by applying $x_*$ to the previously mentioned part of the diagram, we conclude that $x_*\Cl_x\circ x_*x^*\cl_q$ is the identity map on $x_*\Lambda$.
\end{proof}
\end{lemma}

Since we will use this notation more often, we make the following definition.

\begin{definition}
Let $\psi:q^?\Lambda\to x_*\Lambda$ denote the morphism $x_*\Cl_{x}\circ \eta_{x}$ which was used in the above lemma.
\end{definition}

The morphism $\cl_q:q^*\Lambda\to q^?\Lambda$ induces a morphism between the triangle $\Sigma_\X$ and a new triangle $\Sigma^?_\X$, where the new triangle is the localisation triangle associated to $q^?\Lambda$ and $(\X,U,Z)$.
We wish to obtain a factorization $\Sigma_\X\to \Sigma_\X^?\to \Sigma_{\Oh_v}$ of the morphism $\Sigma_\X\to \Sigma_{\Oh_v}$.

\begin{lemma}\label{factorizationtrianglevia?}
There exists a morphism $\varphi:Rj_*j^*q^?\Lambda\to x_*Rj'_*(j')^*\Lambda$ such that the following diagram commutes: 
\begin{equation}\label{openimmersionclasses}
\begin{tikzcd}
q^?\Lambda \arrow[d, "\psi"] \arrow[r, "\eta_j"] & Rj_*j^*q^?\Lambda \arrow[d, "\varphi", dashed] \\
x_*\Lambda \arrow[r, "x_*\eta_{j'}"]       & x_*Rj'_*(j')^*\Lambda            
\end{tikzcd}
\end{equation}
The above commuting diagram extends uniquely to a morphism of triangles $\Sigma_{\X}^?\to \Sigma_{\Oh_v}$.
\begin{proof}
Assume for a moment that one has a commuting diagram as above.
Lemma \ref{uniquetrianglemorphism} shows that the commuting diagram extends to a morphism of triangles and the last statement of Lemma \ref{uniquetrianglemorphism} shows that the morphism of triangles is unique.\par
Now we will construct our morphism $\varphi$. The following diagram of endofunctors on $D^+(\X\et,\Lambda)$ commutes because both ways around equal the adjunction morphism of $j\circ x=x\circ j'$:
\begin{equation}\label{localisationcommutingopen}
\begin{tikzcd}
\Id \arrow[d, "\eta_x"] \arrow[r, "\eta_j"]  & Rj_*j^* \arrow[rd, "Rj_*\eta_xj^*"]                       &               \\
x_*x^* \arrow[r, "x_*\eta_{j'}x^*"] & x_*Rj'_*(j')^*x^* \arrow[r, phantom] \arrow[r, "="] & Rj_*x_*x^*j^*
\end{tikzcd}
\end{equation}

We now set $\varphi$ to be equal to the following composition:
\[Rj_*j^*q^?\Lambda\overset{Rj_*\eta_xj^*}{\to}x_*Rj'_*(j')^*x^*q^?\Lambda\overset{x_*Rj'_*(j')^*\Cl_x}{\to} x_*Rj'_*(j')^*x^?q^?= x_*Rj'_*(j')^*\] It follows easily from commuting Diagram (\ref{localisationcommutingopen}) combined with the description of $\psi$, that $\varphi$ fits into a commuting Diagram (\ref{openimmersionclasses}).
\end{proof}
\end{lemma}

The lemma has the following immediate consequence.

\begin{corollary}\label{Clqandpsi}
The morphism $\alpha:\iota_*\iota^!q^*\Lambda\to x_*\overline{\iota}_*\overline{\iota}^!\Lambda$ factors via $\iota_*\iota^!\cl_q:\iota_*\iota^!q^*\Lambda\to \iota_*\iota^!q^?\Lambda$. The morphism in the factorization $\overline{\psi}:\iota_*\iota^!q^?\Lambda\to x_*\overline{\iota}_*\overline{\iota}^!\Lambda$ is uniquely determined by the commuting diagram:
\begin{equation}\label{psidetermined}
\begin{tikzcd}
\iota_*\iota^!q^?\Lambda \arrow[d, "\overline{\psi}"] \arrow[r, "\epsilon_{\iota}^!"]              & q^?\Lambda \arrow[d, "\psi"] \\
x_*\overline{\iota}_*\overline{\iota}^!\Lambda \arrow[r, "x_*\epsilon_{\overline{\iota}}^!"] & x_*\Lambda               
\end{tikzcd}
\end{equation}
\begin{proof}
One has by Lemma \ref{factorizationtrianglevia?} a factorization $\Sigma_\X\to \Sigma_\X^?\to \Sigma_{\Oh_v}$ of morphisms of triangles, where the first morphism was induced by $\cl_q$. The morphism of triangles $\Sigma_\X^?\to \Sigma_{\Oh_v}$ is unique by Lemma \ref{factorizationtrianglevia?} and so in particular the component $\overline{\psi}:\iota_*\iota^!q^?\Lambda\to x_*\overline{\iota}_*\overline{\iota}^!\Lambda$ is uniquely determined by commuting Diagram (\ref{psidetermined}).
\end{proof}
\end{corollary}

Before we prove the main result we need an auxiliary lemma on compatibility of cycle classes with base-change. We assume that it is well-known, but we could not find a proof in the literature.

\begin{lemma}\label{pullingbackcycleclass}
Let $\varphi:V\to S$ be a finite type morphism of Noetherian $\ZZ[1/n]$-schemes and let $W\subset S$ be a local complete intersection of codimension $c$. Write $\varphi^*:\Ho^{2c}_{|W|}(S,\Lambda(c))\to \Ho^{2c}_{|W\times_S V|}(T,\Lambda(c))$ for the pullback map. Assume that $V\times_SW$ is a local complete intersection of codimension $c$ in $V$. Then we have the equality $\varphi^*\cl(W)=\cl(V\times_S W)$.
\begin{proof}
We first go back to how the cohomology class associated to a codimension $c$ local complete intersection $\alpha:W\hookrightarrow S$ is defined. This is explained in \cite[p.140]{SGA4.5} and we follow this explanation. A priori the cycle class will live in $\Het^0(S,R^{2c}\alpha^!\Lambda(c))$. This last term turns out to be $\Ho^{2c}_{|W|}(S,\Lambda(c))$ by the semi-purity theorem (see for example \cite[Lemma 5.4]{FailureWA}) combined with a spectral sequence.\par 
We explain how to produce the element of $\Het^0(S,R^{2c}\alpha^!\Lambda(c))$. Since the sheaf $R^{2c}\alpha^!\Lambda(c)$ is the sheafification of the presheaf $T\mapsto \Ho^{2c}_{T\times_S W}(T,\Lambda(c))$, it suffices to give elements $s_i\in \Ho^{2c}_{S_i\cap W}(S_i,\Lambda(c))$, where $(S_i)_i$ is an open cover of $S$, such that the $s_i$ agree on overlaps. The agreeing on overlaps happens by construction, therefore we may shrink $S$ so that $W$ is cut out by $c$ equations $f_1,...,f_c\in \Oh_S(S)$ and we are asked to construct an element of $\Ho^{2c}_{|W|}(S,\Lambda(c))$.\par 
Let $D_1,...,D_c$ denote the closed subschemes of $S$ given by $V(f_1),...,V(f_c)$ respectively. Denote the complement of $D_i$ by $U_i$. In $\Het^0(U_i,\Gm)$ we have the element $f_i$. From the localisation exact sequence there is a connecting homomorphism $\delta:\Het^0(U_i,\Gm)\to \Ho^1_{|D_i|}(S,\Gm)$ and associated to the Kummer sequence there is another connecting homomorphism $\delta':\Ho^1_{|D_i|}(S,\Gm)\to \Ho^2_{|D_i|}(S,\Lambda(1))$. We set $\cl(D_i)$ to be $\delta'(\delta(f_i))$. By using the ideas in \cite[Cycle 1.2]{SGA4.5}, we see that there is a cup product $\Ho_{|D_1|}^*(S)\otimes...\otimes \Ho_{|D_m|}^*(S)\to \Ho^*_{|D_1\cap...\cap D_m|}(S)$. Thus we get a section $\cl(W)\in \Ho^{2c}_{|W|}(S,\Lambda(c))$ by setting $\cl(W):=\cl(D_1)\cup...\cup \cl(D_c)$. \par
\hfill \break
We use the above description of the cohomology class of a local complete intersection to prove the lemma. Since the construction of the cycle classes is local, we may already assume that we are in the situation at the end of the third paragraph of this proof, so we have already shrunk $S$. Let $V_i$ denote the open subscheme $V\setminus (D_i\cap V)$ of $V$. The pullback homomorphisms $\varphi^*:\Ho^0(U_i,\Gm)\to \Ho^0(V_i,\Gm)$ sends $f_i$ to $\varphi^*(f_i)$.\par 
The pullback homomorphism $\varphi^*$ is compatible with connecting homomorphisms (this follows from the fact that these connecting homomorphisms come from distinguished triangles and pullback is induced by $\eta_{\varphi}:\Id\to R\varphi_*\varphi^*$). Moreover $\varphi^*$ is compatible with cup products (because the cup product is induced by the tensor product $\Gamma_{|D_i|}(S,\mathcal{F})\otimes \Gamma_{|D_j|}(S,\mathcal{G})\to \Gamma_{|D_i\cap D_j|}(S,\mathcal{F}\otimes \mathcal{G})$ and the tensor product is compatible with $\varphi^*$. These observations lead us to the equality $\varphi^*\cl(W)=\varphi^*\cl(D_1)\cup...\cup \varphi^*\cl(D_c)$, where we have $\varphi^*\cl(D_1)=\delta'(\delta(\varphi^*(f_i))$.\par
By the assumption in the lemma we have that $W\times_S V$ is a local complete intersection of codimension $c$ in $W$, so one can make sense of the cohomology class $\cl(W\times_S V)$. Moreover it is clear that $W\times_S V$ is cut out by the equations $\varphi^*(f_1),...,\varphi^*(f_c)$ in $\Oh_V(V)$. By using the end of the previous paragraph we see that $\cl(W\times_S V)=\varphi^*(\cl(W))$ holds.
\end{proof}
\end{lemma}

By using the lemma we can prove the following statement.

\begin{prop}\label{Psiinsteadofpsi}
The commuting Diagram (\ref{psidetermined}) is induced by a commuting diagram of functors:
\begin{equation}
\begin{tikzcd}
                                                                     & \iota_*\iota^! \arrow[r, "\epsilon_{\iota}^!"] \arrow[ld,  "\overline{\Psi}"', dashed]                     & \Id \arrow[d, "\Psi"] \\
{\iota_*x_{Z,*}x_Z^?\iota^!} \arrow[r, "="] & x_*\overline{\iota}_*\overline{\iota}^!x^? \arrow[r, "x_*\epsilon_{\overline{\iota}}^!x^?"] & x_*x^?               
\end{tikzcd}
\end{equation}
In this diagram the map $\overline{\Psi}$ is depends only on $\cl(B)\in \Ho^{2d}_{|B|}(Z,\Lambda(d))$.
\begin{proof}
All of the functors in this diagram are representable by Lemma \ref{representHZ}. So it will suffice to show that there is a commuting diagram as follows:
\begin{equation}\label{representedbyLambda}
\begin{tikzcd}
\Lambda_Z                                         & \Lambda \arrow[l, "\eta_{\iota}"]                                                \\
{\Lambda_B(-d)[-2d]} \arrow[u, "\overline{\Phi}", dashed] & {\Lambda_{x}(-d)[-2d]} \arrow[u, "\Phi"] \arrow[l, "\eta_{\overline{\iota}}"]
\end{tikzcd}
\end{equation}
Here $\Phi$ is the homomorphism corresponding to $\cl_x$ and with $\Lambda_{x}$ we mean $\Lambda_{x(\Sp(\Oh_v))}$. By fully faithfulness of the functor $\iota_*$ we have $\Hom_{D^+(\X)}(\Lambda_{B}(-d)[-2d],\Lambda_Z)=\Hom_{D^+(Z)}(\Lambda_B,\Lambda_Z(d)[2d])=\Ho^{2d}_{|B|}(Z,\Lambda(d))$, where the last equality comes from using adjunction from $x_Z$ on $Z\et$. Any morphism $\Lambda_{x}(-d)[-2d]\to \Lambda_Z$ factors uniquely over $\Lambda_{x}(-d)[-2d]\to \Lambda_B(-d)[-2d]$ (look at the stalks), so we obtain an identification $\Hom_{D^+(\X)}(\Lambda_{x}(-d)[-2d],\Lambda_Z)=\Ho^{2d}_{|B|}(Z,\Lambda(d))$.\par
If we go the other way around the diagram, we see that $\eta_{\iota}$ corresponds to pullback in cohomology. The map $\Phi:\Lambda_{x}(-d)[-2d]\to \Lambda$ corresponds to the element $\cl(x)\in \Ho^{2d}_{|x|}(\X,\Lambda(d))$. The element $\eta_\iota\circ \Phi $ corresponds to $\iota^*(\cl(x))$ and by Lemma \ref{pullingbackcycleclass} we get the equality $\iota^*(\cl(x))=\cl(B)\in \Ho^{2d}_{|B|}(Z,\Lambda(d))$.\par
So we see that by defining $\overline{\Phi}$ as the morphism that corresponds under adjunction to $\cl(B)$ makes Diagram~(\ref{representedbyLambda}) commute. By representability this gives us the commuting diagram that was claimed in the proposition.
\end{proof}
\end{prop}

Now we can conclude the proof of main theorem.

\begin{theorem}\label{MainTH}
Consider the setting of Diagram (\ref{setting}). The evaluation map $u^*:\br(U)[n]\to \br(k_v)[n]$ depends only on $\cl(Z\cap x)\in \Ho^{2d}_{|Z\cap x|}(Z,\mu_n^{\otimes d})$.
\begin{proof}
By Lemma \ref{determinedbyalpha}, we have that $u^*:\Het^2(U,\mu_n)\to \Het^2(k_v,\mu_n)$ is determined by $\alpha:\iota_*\iota^!q^*\Lambda\to x_*\overline{\iota}_*\overline{\iota}^!\Lambda$. By Corollary \ref{Clqandpsi}, we see that $\alpha$ is determined by $\cl_q$ and the map $\overline{\psi}:\iota_*\iota^!q^?\Lambda\to x_*\overline{\iota}_*\overline{\iota}^!\Lambda$. By Proposition~\ref{Psiinsteadofpsi} we have that $\overline{\psi}$ is determined by $\cl(B)\in \Ho^{2d}_{|B|}(Z,\Lambda(d))$, where we have $B=Z\cap x$. 
\end{proof}
\end{theorem}

We show why Theorem \ref{snc theorem} implies a special case of Theorem \ref{MainTH}.

\begin{remark}\label{snctheoremcomparison}
Assume that $Z$ is an snc divisor and that we have $\gcd(n,q-1)=1$, which is the setting of Theorem \ref{snc theorem}. Denote the irreducible components of $Z$ by $(Z_i)_i$. Suppose that $Z_i$ meets $x$ with multiplicity~$m_i$. Gabber's purity theorem \cite{TravauxdeGabber}[Exposé XVI.3 Théorème 3.1.1] implies that for all $i$ with $m_i\geq 1$, we have a canonical isomorphism $\ZZ/n\to \Ho^{2d}_{|B|}(Z_i,\Lambda(d))$ and the cycle class of $x\cap Z_i$ is mapped to $m_i$. In particular we see that under the composition 
\begin{equation}\label{sncrestriction}\Ho^{2d}_{|B|}(Z,\Lambda(d))\to \bigoplus_{i,\,m_i\geq 1}\Ho^{2d}_{|B|}(Z_i,\Lambda(d))\overset{\sim}{\leftarrow} \bigoplus_{i,\,m_i\geq 1}\ZZ/n\end{equation}
the element $\cl(B)$ is mapped to $(m_i)_i$. So we see that by Theorem \ref{snc theorem}, the evaluation map $u^*$ is determined by $\cl(B)$, since it is even determined by the coarser invariants $(m_i)_i$ and the reduction of $x$ modulo $v$.
\end{remark}

The map (\ref{sncrestriction}) need not be injective as shown in this example.

\begin{example}
We return to Example \ref{babyexample} one last time, so we have $Z=\Sp(\ZZ_p[t]/(tp))$ with components~$\A^1_{\F_p}$ and $\Sp(\ZZ_p)$ and in this situation $|B|$ corresponds to the maximal ideal $(p,t)$. We wish to show that whenever $\gcd(n,p-1)\neq 1$, the map $\Ho^2_{|B|}(Z,\mu_n)\to \Ho^2_{|B|}(\A^1_{\F_p},\mu_n)\oplus \Ho^2_{|B|}(\ZZ_p,\mu_n)=(\ZZ/n)^2$ is not injective. We use the localisation exact sequence for $(Z,Z\setminus |B|)$. By easy calculations we get \[\Het^1(Z,\mu_n)=\Oh(Z)^\times/(\Oh(Z)^\times)^n\cong \F_p^\times/(\F_p^\times)^n\] \[\Het^1(Z\setminus |B|,\mu_n)=k_v^\times/(k_v^\times)^n\oplus (\F_p[t,1/t])^\times/((\F_p[t,1/t])^\times)^n,\] which is isomorphic to two copies of $\ZZ/n\oplus \F_p^\times/(\F_p^\times)^n$. The inclusion of $\Het^1(Z,\mu_n)$ into $\Het^1(Z\setminus |B|,\mu_n)$ is identified with the diagonal inclusion of $\F_p^\times/(\F_p^\times)^n$ into its self-product (and the inclusion is zero on the $(\ZZ/n)^2$-component). It follows that the kernel of $\Ho^2_{|B|}(Z,\mu_n)\to \Ho^2_{|B|}(\A^1_{\F_p},\mu_n)\oplus \Ho^2_{|B|}(\ZZ_p,\mu_n)$ has a subgroup isomorphic to $\F_p^\times/(\F_p^\times)^n$, which is nonzero by our assumption on $n$.
\end{example}

This motivates us to ask the following question of which we do not know the answer.

\begin{question}In the setting of Remark \ref{snctheoremcomparison} does the condition $\gcd(n,q-1)=1$ guarantee that the map $\Ho^{2d}_{|B|}(Z,\Lambda(d))\to \bigoplus_{i,\,m_i\geq 1}\Ho^{2d}_{|B|}(Z_i,\Lambda(d))$ is an isomorphism?
\end{question}

We close with a remark that shows that a certain converse of Theorem \ref{slogan} is not true.

\begin{remark}
Suppose that the evaluation maps $u_1^*$ and $u_2^*$ on $\br(U)$ agree.
Do we then have $|x_1\cap Z|=|x_2\cap Z|$ as well as the equality $\cl(x_1\cap Z)=\cl(x_2\cap Z)$ in $\Ho^{2d}_{|x_1\cap Z|}(Z,\mu_n(d))$ for all $n$? \par This statement is false by an almost trivial counterexample. Namely take $U=\An^1_{\Q_p}\overset{\pi}{\to}\Sp(\Q_p)$ sitting inside $\X=\Pn^1_{\ZZ_p}$ with the boundary being $Z:=\Pn^1_{\F_p}\cup \{\infty\}_{\ZZ_p}$ in this case. The natural map $\pi^*:\br(\Q_p)\to \br(\An^1_{\Q_p})$ is an isomorphism. So any point $u\in \An^1(\Q_p)$ has that the evaluation map $u^*$ equals the inverse of $\pi^*$. However, for $x$ the unique lift of $u$ we do not have that $|x\cap Z|$ is the same for every point $u\in \An^1({\Q_p})$.
\end{remark}

\printbibliography

\end{document}